%% file: circdiff.tex
\documentstyle[fullpage, 12pt]{article}
\input{imsmark.tex}
\input{psfig}
\begin{document}
\title{Topological Conjugacy of Circle Diffeomorphisms}
\author{JUN HU\\Department of Mathematics, Graduate Center \\CUNY,
NY10036, USA \\{} \\DENNIS P. SULLIVAN\\Department of Mathematics,
Graduate Center \\CUNY, NY10036, USA and IHES, FRANCE}
\maketitle
\SBIMSMark{1995-7}{May 1995}{}
\begin{abstract}
The classical criterion for a circle diffeomorphism to be topologically
conjugate to an irrational rigid rotation was given by A. Denjoy [1].
In [5] one of us gave a new criterion. There is an
example satisfying Denjoy's bounded variation condition rather than
[5]'s Zygmund condition and vice versa. This paper will give
the third criterion which is implied by either of the above criteria.
\end{abstract}

\thispagestyle{empty}
\centerline{\bf Contents}
\begin{enumerate}
\item Introduction
\item Cross ratio distortion
\item Non wandering set and ergodicity
\item Proof of results
\item Three examples
\item Appendix
\end{enumerate}
\section{Introduction}
Given a circle orientation preserving homeomorphism $f:
S^{1}\rightarrow S^{1}$,
the rotation
number \[\rho(f)=\lim_{n\rightarrow \infty}\frac{F^{n}(x)-x}{n} \;\;mod
\;\; 1\]
is independent of $x$ and the lift $F$ of $f$, where
$F:R^{1}\rightarrow R^{1}$ is a lift of $f$ and $x\in
R^{1}$. And it is invariant under
topological conjugations. The rotation number $\rho(f)$ is
a rational number if and only if $f$ has a periodic orbit.
From the theory of Poincar\'e, for an orientation preserving homeomorphism
$f:S^{1}\rightarrow S^{1}$, if $f$ has
a periodic orbit then its dynamics turn out trivial: any two periodic
orbits have the same period and any orbit tends to a periodic orbit; if
$f$ doesn't have any periodic orbit then it is semi-conjugate to an
irrational rigid
rotation. A natural question is whether or not the semi-conjugation
could be improved to be a topological conjugation. In the
following context when we say a rigid rotation we always mean an
irrational rigid rotation. Denjoy proved the following.
\begin{description}\item[Theorem A] Given an orientation preserving
homeomorphism $f$ of the circle $S^{1}$ with an irrational rotation
number,
$f$ is topologically conjugate to a rigid rotation provided $f$ is
$C^{1}$ and the logarithm of the
derivative of $f$ is of bounded variation.
\end{description}
There is an example (Denjoy
counterexample) to show that $C^{1}$ smoothness is not enough [6].
It is shown even $C^{\infty }$ smoothness is not enough yet in [15].
Actually
an orientation preserving circle homeomorphism with an irrational
rotation
number is topologically conjugate to an irrational rigid rotation if
and only if it has no wandering interval. Denjoy achieved this by
controlling the variation of the derivative.
Recently one of us proved the non existence of wandering interval by
assuming the logarithm of the derivative satisfies the Zygmund
condition.
\begin{description}\item[Definition]
A continuous map $f:R^{1}\rightarrow R^{1}$ satisfies the Zygmund
condition if there
exists $B>0$ such that \[\sup_{x,t}|\frac{f(x+t)+f(x-t)-2f(x)}{t}| \leq
B.\]
\end{description}
\begin{description}\item[Theorem B] Given an orientation preserving
homeomorphism $f$ of the circle $S^{1}$ with an irrational rotation
number, $f$
is topologically conjugate to a rigid rotation if $f$ is $C^{1}$ and the
logarithm of the derivative satisfies the Zygmund condition.
\end{description}
But there is an example satisfying Denjoy's bounded variation condition
and not Zygmund's condition and vice versa [section 5]. This paper
gives a third criterion which is implied by either of the above
two and which implies $f$ is topologically conjugate to a rigid
rotation.
\begin{description}\item[Definition] Let $I$ be a closed interval of
$R^{1}$.
A continuous map $f:I\rightarrow R^{1}$ is of bounded Zygmund variation
if there exists $B>0$ such that
\[\sup_{\{x_{0},x_{1},\cdots,x_{n}\}}\sum_{i=0}^{n-1}
|f(x_{i})+f(x_{i+1})-2f(\frac{x_{i}+x_{i+1}}{2})| \leq B ,\]
where $\{ x_{0},x_{1},\cdots,x_{n} \} $ is a partition of the interval
$I$. The supremum is called Zygmund variation of $f$ over $I$. It is denoted by
$ZV(f|_{I})$.
\end{description}
\begin{description}\item[Definition] Let $I$ be a closed interval of
$R^{1}$. A continuous map $f:I\rightarrow R^{1}$ is of bounded quadratic
variation if there exists $B>0$ such that
\[\sup_{\{x_{0},x_{1},\cdots,x_{n}\}
}\sum_{i=0}^{n-1}(f(x_{i+1})-f(x_{i}))^{2} \leq B ,\]
where $\{x_{0},x_{1},\cdots,x_{n}\}$ is a partition of the interval
$I$. The supremum is called quadratic variation of $f$ over $I$. It is denoted by $QV(f|_{I})$.
\end{description}
\begin{description}\item[Theorem C] Given an orientation preserving
homeomorphism $f$ of the circle $S^{1}$ with an irrational rotation
number, $f$
is topologically conjugate to a rigid rotation if $f$ is $C^{1}$ and the
logarithm of the derivative has bounded Zygmund variation and
bounded quadratic variation.
\end{description}
\section{Cross ratio distortion}
In this section we control cross ratio distortion for
standard 4-tuples in terms of Zygmund variation and quadratic variation
(compare $\S 1$ of [5]).
Let $a,b,c,d\in R^{1}$ and $a<b<c<d$.

One cross ratio $[a,b,c,d]=\frac{(d-b)(c-a)}{(c-b)(d-a)}$ can be
computed by \[log[a,b,c,d]=\int \int_{S}\frac{dxdy}{(x-y)^{2}},\]
where $S$ is $\{(x,y):a\leq x\leq b, c\leq y\leq d\}$.

Another cross ratio $(a,b,c,d)$ is $\frac{(c-b)(d-a)}{(b-a)(d-c)}$
and, obviously, \[[a,b,c,d]=1+\frac{1}{(a,b,c,d)}.\]

Given a homeomorphism $h$, the distortion of the second cross ratio
under $h$ is \[\frac{(ha,hb,hc,hd)}{(a,b,c,d)}.\]

In this paper, by the cross ratio distortion we mean the distortion of
the second cross ratio.

We call a 4-tuple $a<b<c<d$ standard if $b-a=c-b=d-c$. The cross ratio
distortion under $h$ of a standard 4-tuple is bounded away from zero
and from above if and
only if $(ha,hb,hc,hd)$ is also. If $h$ is $C^{1}$ diffeomorphism, then
\[log[1+\frac{1}{(ha,hb,hc,hd)}]=log[ha,hb,hc,hd]=\int \int _{S}(h\times
h)^{*}\mu ,\]
where $\mu $ is the measure $\frac{dxdy}{(x-y)^{2}}$.

Clearly the cross ratio distortion under $f$ of a standard 4-tuple is
bounded
away from zero and from above if and only if $log[ha,hb,hc,hd]$ is also.
Calculating the integrand, we get \[\frac{h'xh'y}{(hx-hy)^{2}}=
\frac{1}{(x-y)^{2}}\frac{h'xh'y}{[h']^{2}_{xy}},\]
where $[h']_{xy}$ is the average of $h'$ over the interval $[x,y]$.

Since $b-a=c-b=d-c$,
$\int \int_{S}\frac{dxdy}{(x-y)^{2}}=log([a,b,c,d])=log\frac{4}{3}$.
Thus a bound on $\frac{h'xh'y}{[h']^{2}_{xy}}$ yields a bound on
the cross ratio distortions for standard 4-tuples.

We say $h$ satisfies the {\em bounded Koebe condition} if one of the
following equivalent conditions hold:
\[1)\;\;\frac{1}{M}\leq \frac{h'xh'y}{[h']^{2}_{xy}} \leq M \; for\;
some\; M>0,\]
\[2)\;\;|log\frac{h'xh'y}{[h']^{2}_{xy}}|\leq M' \; for\; some
\;M'>0.\]
The following proposition is trivial.
\newtheorem{prop}{Prop.}
\begin{prop}
If $h$ satisfies the bounded Koebe condition then the cross ratio
distortion under $h$ of a standard 4-tuple is bounded away
from zero and from above.
\end{prop}
In order to estimate the $\log $ in 2), i.e.,
\[logh'x+logh'y-2log[h']_{xy},\]
let us consider the following two terms:
\[a)\;\;logh'x+logh'y-2[logh']_{xy}\] and
\[b)\;\;log[h']_{xy}-[logh']_{xy}\;.\]
Remark: If both a) and b) are bounded, then 1) and 2) hold.

Expression a) can be controlled by the Zygmund
variation of $logh'$ on the interval $[x,y]$ because of the following
proposition.
\begin{prop}
Let $\phi $ be a continuous function from $R^{1}$ to $R^{1}$. Then
\[|\phi (x)+\phi (y)-2[\phi ]_{xy}|\]
is no more than the Zygmund variation $ZV(\phi |_{[x,y]})$ of $\phi corollary$ over $[x,y]$.
\end{prop}
Remark: As we define the Zygmund variation of $\phi $ on the interval 
$[a, b]$ in the introduction, we can also define the average Zygmund 
variation of $\phi $ on $[a, b]$ by replacing the value $\phi (\frac{x_{i}+
x_{i+1}}{2})$ of $\phi 
$ at the middle point by the average $\frac{1}{|x_{i+1}-x_{i}|}
\int _{[x_{i}, x_{i+1}]}\phi $ of $\phi $ over $[x_{i}, x_{1+1}]$. 
Prop. 2 tells us that the average Zygmund variation of $\phi $ over 
$[a, b]$ is no more than the Zygmund variation of $\phi $ over $[a, b]$.
Conversely one can show that the Zygmund variation of $\phi $ over $[a, b]$
is no more than twice the average Zygmund variation of $\phi $ over
$[a, b]$.  
Hence these two conditions are actually equivalent.

Proof: Without loss of generality, assume $[x,y]=[0,1]$.
Then \[[\phi ]_{01}=\int _{0}^{1}\phi dx. \]
Suppose that we define successive approximations to the average of $\phi $
over $[a,b]$ by
\[A_{0}[a,b] = (\phi (a) + \phi (b))/2\]
and
\[A_{n+1}[a,b] = (A_{n}[a,m] + A_{n}[m,b])/2 ,\]
where  $m = (a+b)/2 $. Similarly, measure non-linearity by expressions
\[N_0[a,b] = A_0[a,b] - A_1[a,b] = (\phi(a) - 2 \phi(m) + \phi(b))/4\]
and 
\[N_{n+1}[a,b] =(N_n[a,m] + N_n[m,b])/2 ,\]
or equivalently
\[N_{n}[a,b] = A_{n}[a,b] - A_{n+1}[a,b] .\]
Then
\[A_{0}[0,1] -A_{n}[0,1] = N_{0}[0,1] + \cdots  + N_{n-1}[0,1]\]
with
\[| N_{k}[0,1] | \leq ZV(\phi|_{[0,1]})/2^{k+2}\]
hence	
\[2 |A_{0}[0,1] - \lim _{n\rightarrow \infty } A_{n}[0,1]| 
        \leq ZV(\phi |_{[0,1]}) ,\]
i.e., 
  \[|\phi (0)+\phi (1)-2[\phi ]_{01}|\leq ZV(\phi |_{[0,1]}).\]

Next we estimate the expression b) in terms of the quadratic variation.
\newtheorem{lemma}{Lemma}
\begin{lemma}
If $\epsilon \geq \delta >-1$, assume \[log(1+\epsilon)=\epsilon
-\frac{\epsilon ^{2}}{2}\Delta (\epsilon ),\]
then there exists $B(\delta )>0$ depending on $\delta $ such that
$|\Delta (\epsilon )|\leq B(\delta )$.
\end{lemma}
Proof: Since \[\Delta (\epsilon )=\frac{\epsilon -log(1+\epsilon )
}{\epsilon ^{2}/2},\]
The proof is an elementary calculation.
\newtheorem{definition}{Definition}
\begin{definition}
A quantity $C_{1}$ (depending on parameters) is a big O of another
quantity $C_{2}$ (depending on the parameters) if there exists a
constant $B$ (independent of the parameters) such that
\[|C_{1}|\leq B|C_{2}|.\]
\end{definition}
\begin{prop}
Suppose the derivative $h'$ satisfies $1/C\leq h' \leq C $ for some
$C>0$. Then
the expression b) is equal to the big $O$ of the quadratic variation
of $logh'$ over the interval $[x,y]$.
\end{prop}
Proof: Let $h'(x)=a$. The expression b) is unchanged if we multiply
$h'x$ by $1/a$. Write $(1/a)h'$ on $J=[x,y]$ as $1+\epsilon $ where
$\epsilon $ is a function of $(t-x), t\in J$. Expand the two terms of b)
\[log\frac{1}{|J|} \int_{J}(1+\epsilon)
-\frac{1}{|J|}\int_{J}log(1+\epsilon )\]
\[=log(1+\frac{1}{|J|} \int_{J}\epsilon )
-\frac{1}{|J|}\int_{J}[\epsilon -\frac{\epsilon ^{2}}{2}\Delta
(\epsilon )]\]
\[=[\frac{1}{|J|}\int_{J} \epsilon -\frac{1}{2}(\frac{1}{|J|}\int_{J}
\epsilon)^{2} \Delta (\frac{1}{|J|}\int_{J} \epsilon
)]-[\frac{1}{|J|}\int_{J} \epsilon -\frac{1}{|J|}\int_{J}\frac{\epsilon
^{2}}{2}\Delta (\epsilon )]\]
\[=-\frac{1}{2}(\frac{1}{|J|}\int_{J}\epsilon )^{2}\Delta
(\frac{1}{|J|}\int_{J}\epsilon )+\frac{1}{|J|}\int_{J}\frac{\epsilon
 ^{2}}{2}\Delta (\epsilon ).\]
By the Cauchy inequality, $(\frac{1}{|J|}\int _{J}\epsilon )^{2}\leq
\frac{1}{|J|}\int _{J}\epsilon ^{2}$.
Since $1/C\leq h'\leq C$ for some $C>0$, there exists $\delta (C)>-1$
such that $\epsilon =\frac{h^{'}t}{h^{'}x}-1\geq \delta $ for any
$t\in J$, $J=[x, y]$.
Hence $\frac{1}{|J|}\int _{J}\epsilon \geq \delta $. By
the Lemma 1, there exists $B(\delta )>0$ such that $|\Delta (\epsilon )|
\leq B(\delta )$. Hence $|\Delta (\frac{1}{|J|}\int _{J}\epsilon )|\leq
B(\delta )$.
Furthermore we can get that $\epsilon =\frac{h't}{h'x}-1$ is
a big $O$ of $log\frac{h't}{h'x}=logh't-logh'x$. so the
expression b) is a big $O$ of the quadratic variation of $logh'$ over
$J$.
The following proposition will be used in section 4 to the
iterates of a circle diffeomorphism $f$.
\begin{prop}
Suppose $h:I\rightarrow R^{1}$ is a $C^{1}$ diffeomorphism with $h'>0$,
and $logh'$ has bounded Zygmund variation and bounded quadratic
variation over $I$. Assume $J_{0}\subset I$ and $J_{0}\;
J_{1}=h(J_{0}),\;\cdots,\;J_{n}=h^{n}(J_{0})$ are pairwise
disjoint. Then the cross ratio distortion under $h^{n}$ of a standard
4-tuple in the interval $J_{0}$ is the big $O$ of the sum of the
Zygmund variation and the quadratic variation of $\log h^{'}$ on
$\cup _{i=0}^{n-1}J_{i}$.
\end{prop}
Proof: From the expression 2) above the Prop. 1, we want to
estimate \[log
\frac{(h^{n})^{'}(x)(h^{n})^{'}(y)}{[(h^{n})^{'}]^{2}_{xy}}.\]
By the chain rule of calculating the derivative of $h^{n}$,
\[log\frac{(h^{n})^{'}(x)(h^{n})^{'}(y)}{[(h^{n})^{'}]^{2}_{xy}}
=\sum _{i=0}^{n-1}log
\frac{h^{'}(h^{i}(x))h^{'}(h^{i}(y))}{[h^{'}]^{2}_{h^{i}(x)h^{i}(y)}}.\]
Each summand can be decomposed into the expression a) and
expression b), by the Prop. 2 and Prop. 3, each summand is the big
$O$ of the sum of the Zygmund variation and the quadratic variation of
$logh^{'}$ over the interval $[h^{i}(x), h^{i}(y)]$, where $i=0, 1, 2,
\cdots , n-1$.
So the cross ratio distortion under $h^{n}$ of a standard 4-tuple in
$J_{0}$ is the big $O$ of the sum of the Zygmund variation and the
quadratic variation of
$logh^{'}$ on $\cup _{i=0}^{n-1}J_{i}$.
\section{Nonwandering set and ergodicity}
In this section we review some basic techniques due to Denjoy [1].
Suppose $f:S^{1}\rightarrow S^{1}$ is an orientation preserving
homeomorphism with an irrational rotation number.

For $x\in S^{1}$, let \[\omega (x)=\cap_{n\in N}Cl(\{f^{k}(x)|k\geq
n\}),\] \[\alpha (x)=\cap_{n\in N}Cl(\{f^{-k}(x)|k\geq n\}),\]
where $Cl(A)$ means the closure of the set $A$.
They are called $\omega $ limit set of the orbit of $x$ and
$\alpha $ limit set of the orbit of $x$ respectively.

$x\in S^{1}$ is called a wandering point of $f$ if there exists a
neighborhood $U$ of $x$ such that \[f^{k}(U)\cap U=\emptyset , \forall
k\in Z\setminus \{0\}.\]
A point is called a nonwandering point if it is not a wandering point.
$\Omega (f)$ denotes the set of all nonwandering points, which is
called nonwandering set. Clearly it is a closed subset.

A subset $A$ is invariant under $f$ if \[f(A)\subset
A,\;f^{-1}(A)\subset A.\]
A non-empty subset $A$ is minimal for $f$ if it is closed, invariant
under $f$ and there is no non-empty proper closed subset of $A$ which is
invariant under $f$.
\begin{prop}
Suppose $f$ has no periodic point, then

(1) $\Omega (f)=\omega (x)=\alpha (x),\;\forall x\in S^{1};$

(2) $\Omega (f)$ is a minimal set of $f$;

(3) either $\Omega (f)$ is a nowhere dense perfect subset of $S^{1}$ or
 $\;\Omega (f)=S^{1}$.
\end{prop}
Proof: (1) $\omega (x)$ is a non-empty closed invariant subset of
$S^{1}$. Let $(\gamma ,\delta )$ be a component of
$S^{1}\setminus \omega
(x)$, then $f^{j}((\gamma ,\delta ))$ is also a component of
$S^{1}\setminus \omega (x)$ for any $j\in Z$. Since $f$ has no periodic
point, $\{ f^{j}([\gamma ,\delta ])| j\in Z\} $ must be pairwise
disjoint and hence $(\gamma,\delta)$ is a wandering interval
of $f$. So $S^{1}\setminus \omega (x) \subset
S^{1}\setminus \Omega (f)$. Hence $\Omega (f)\subset \omega (x)$.
Clearly $\omega (x)\subset \Omega (f)$. So $\Omega (f)=\omega (x)$.
Similarily $\Omega (f)=\alpha (x)$.

(2) Clearly from (1).

(3) Let $\partial \Omega $ denote the boundary of $\Omega $, $\partial
\Omega $ is closed. Since \[\partial \Omega \subset \Omega ,\;f(\partial
\Omega )=\partial f(\Omega )=\partial \Omega ,\]
either $\partial \Omega =\emptyset $ hence $\Omega =S^{1}$ or $\partial
\Omega =\Omega $ hence $\Omega $ is nowhere dense. For the second case,
 $\Omega $ is perfect since $\Omega =\omega (y),\; \forall y\in \Omega
$.
\begin{definition}
Suppose $f$ has no periodic point. We say $f$ is ergodic if $\Omega
(f)=S^{1}$, otherwise we say $f$ is not ergodic.
\end{definition}
The following result is well known and its proof can be found in
several references ([1], [2], [3], [4] and etc.).
\begin{prop}
Suppose an orientation preserving homeomorphism $f:S^{1}\rightarrow
S^{1}$ has no periodic point and is ergodic, $\alpha =\rho (f)$. Then
$f$ is topologically conjugate to an irrational rigid rotation
$\tau _{\alpha }:S^{1}\rightarrow S^{1}$ given by \[\tau _{\alpha }(\xi
) =e^{2\pi i\alpha }\xi .\]
\end{prop}
\section{Proofs of results}

A circle homeomorphism with an irrational rotation number is
topologically conjugate to a rigid rotation if and only if it is
ergodic, in other words if and only if it has no wandering interval.
Denjoy's $C^{1+b.v}$-condition and [5]'s $C^{1+Z}$-condition both
guarantee the nonexistence of a wandering interval. In this section we
prove that the $C^{1}$-plus bounded Zygmund variation and bounded
quadratic variation guarantee the nonexistence of a wandering interval.
Before we get into the proofs of these results, we need the following
technique lemmas.
\begin{prop} [Contraction Principle]([10], [6])
Suppose $f: S^{1}\rightarrow S^{1}$ is a circle
homeomorphism has no periodic orbits and $I$ is a subinterval of
$S^{1}$. If $\inf _{n\geq 0} \{|f^{n}(I)|\}=0$, then $I$ is a wandering
interval of $f$.
\end{prop}
Proof: Let $I_{n}= f^{n}(int I)$ and $\Sigma =\cup _{n\geq 0}I_{n}$.

Case 1: If $\Sigma =S^{1}$, then $S^{1}$ is covered by finite
$I_{n_{i}}, i=1, 2, ..., k$. Since $\inf _{n\geq 0}{|I_{n}|}=0$,
the Lebesgue's lemma implies that there is $I_{l}$, $l\in N$, contained
in one of $I_{n_{i}}, i=1, 2, ..., k$, name it $I_{j}$. Without loss of
generality we assume $l>j$. Note $I_{l}=f^{l-j}(I_{j})$. So $f^{l-j}$
will have a periodic point in the closure of $I_{j}$. This is
contradiction.

Case 2: Suppose $\Sigma \neq S^{1}$. If there is no component $U$ of
$\Sigma $ such that some iterate of $U$ intersects with $U$, then $U$
and hence $I$ are wandering intervals. If there is a component $U$ of
$\Sigma $ such that $f^{n}(U)\cap U\neq \emptyset $ for some $n\geq 0$,
then $f^{n}(U)\subseteq U$ and hence $f^{n}$ has a periodic point in the
closure of $U$. This is again a contradiction.

\begin{lemma}[Real Koebe Principle]
If $h:I\rightarrow R^{1}$ does not increase the cross ratio distortions
for standard 4-tuples too much then the quasisymmetric distortions
for standard interior triples are controlled.
More precisely, if $x,\; y\in I$ satisfy $|x-y|$ is as small as the
distance to the boundary $\partial I $ of $I$ and $z=(x+y)/2 $, then
\[\frac{1}{C}\leq |h(x)-h(z)|/|h(z)-h(y)| \leq C,\]
where $C$ only depends on the bound of the cross ratio distortions for
standard 4-tuples.
\end{lemma}
Proof: See $\S 2$ of [5]. The idea to prove this lemma is to use the
four interval arguement.
Let $J, L, M, R$ be four contiguous equal lenth
intervals. Suppose the lenth of $h(L)$ is much smaller than $h(M)$.
Since the cross ratio distortion
$\frac{|h(M)||h(T)|}{h(L)||h(R)|}/3$ on $L, M, R$ is greater than the
ratio distortion $\frac{|h(M)|}{|h(L)|}$, no bound of ratio distortions
implies no bound of cross ratio distortions.

\vspace{.2in}
It is easy to use the Real Koebe Principle to get the following
Macroscopic Koebe Distortion Principle.
\begin{definition}
Let $M$ and $T$ be two intervals with $M\subset T$, and $L$ and $R$
be components of $T\setminus M$. If $\epsilon >0$
we say $T$ is an $\epsilon $-scaled
neighborhood of $M$ if
\[\frac{|L|}{|M|}\geq \epsilon \;and \;\frac{|R|}{|M|}\geq \epsilon .\]

\end{definition}
\begin{prop} [Macroscopic Koebe Distortion Principle]

Given any $B>0$, $\epsilon >0$, there exists $\delta >0$ only
depending on $B$ and $\epsilon $ such that,
for any homeomorphism $f$ of the circle, any subintervals $M\subset T$
and any $n\geq 0$, if the cross ratio distortion under $f^{n}$ of any
standard 4-tuple in $T$ is bounded by $B$ and $f^{n}(T)$ contains an
$\epsilon $-scaled neighborhood of $f^{n}(M)$ then $T$ contains a
$\delta $-scaled neighborhood of $M$.

\end{prop}
Proof: Let $T\setminus M=L\cup R$. Without loss of generality, we only
need to prove $\frac{|M|}{|L|}$ can not be very large. Suppose
$\frac{|M|}{|L|}$ is large, we cut $M$ into pieces  $L_{i}$ from left to
right with lengths $2^{i-1}|L|, i=1, 2, 3, \cdots $. We also denote
$L_{0}=L$. From the Real Koebe Principle, there exists a constant $C$
only depending $B$ such that
\[\frac{|f^{n}(L_{i})|}{|\cup _{j=0}^{i-1}f^{n}(L_{j})|}\geq
\frac{1}{C},\]
where $i=1, 2, 3, \cdots $. Hence
\[\frac{|\cup _{j=0}^{i}f^{n}(L_{i})|}
{|\cup _{j=0}^{i-1}f^{n}(L_{j})|}\geq
1+\frac{1}{C},\]
where $i=1, 2, 3, \cdots $. So
\[\frac{|\cup _{j=0}^{i}f^{n}(L_{i})|}
{|f^{n}(L_{0})|}\geq
(1+\frac{1}{C})^{i},\]
where $i=1, 2, 3, \cdots $. This means
\[\frac{|\cup _{j=1}^{i}f^{n}(L_{i})|}
{|f^{n}(L_{0})|}\geq
(1+\frac{1}{C})^{i}-1,\]
where $i=1, 2, 3, \cdots $.

Clearly $i$ can not be very large, otherwise $f^{n}(T)$ can not be an
$\epsilon $-scaled neighborhood of $f^{n}(M)$. Hence we can find a bound
of $i$ only depending on $B$ and $\epsilon $, which means there exists
$\delta >0$ only depending on $B$ and $\epsilon $ such that $T$ contains
a $\delta $-scaled neighborhood of $M$.
\begin{definition}
The intersection multiplicity of a collection of sets
$X_{\alpha \in \Lambda }$ is the maximal cardinality of a
subcollection with non-empty intersection.
\end{definition}
Use the Contraction Principle (Prop. 7), it is easy
to get the following proposition.
\begin{prop}
Suppose $f:S^{1}\rightarrow S^{1}$ is an orientation preserving
homeomorphism without
periodic orbits. Let $I$ be a wandering interval and
not contained in any larger wandering interval. If $I$ is a proper subset
of an interval $J$, then the intersection
multiplicity of the pullbacks $\{f^{-i}(J): i=0, 1, 2, ...\}$ is
infinity.
\end{prop}
Proof: Suppose the intersection multiplicity of the pullbacks 
$\{f^{-i}(J): i=0, 1, 2, ...\}$ is finite. Then
$|f^{-n}(J)|\rightarrow 0$ as $n\rightarrow \infty $. Now apply the 
contraction Principle to $J$ and the map $f^{-1}$. It says that $J$ is a wandering interval. But this is false because $I$ was a maximal wandering interval.

\begin{definition}
Let $f:S^{1}\rightarrow S^{1}$ be an orientation preserving
homeomorphism. The variation of the logarithm of cross ratio distortion
under $f$ is defined as
\[\sup _{\{x_{0},x_{1},\cdots ,x_{n}\}}\sum _{i=0}^{n-1}\sup
_{b_{i},c_{i}\in (x_{i},x_{i+1})}\log \frac{(f(x_{i}), f(b_{i}), f(c_{i}), f(x_{i+1}))}{(x_{i},b_{i},c_{i},x_{i+1})}\]
where $b_{i}$ and $c_{i}$ belong to the open interval
$(x_{i},x_{i+1})$ from $x_{i}$ to
$x_{i+1}$ counter clockwisely and $\{x_{0},x_{1},\cdots ,x_{n}\}$ is a
partition of $S^{1}$.
\end{definition}

Let $f:S^{1}\rightarrow S^{1}$ be an orientation preserving
homeomorphism with an irrational rotational number. Let $I$ be a
wandering interval for $f$. The following combinatorial machinery on 
wandering intervals, $I_{n}=f^{n}(I): n=0,1,2,...$, was
developed in [10] and can be found in [6]. 
\begin{definition} If $n\in N$, we say $I_{k}$ is a left (or right)
predecessor
of $I_{n}$ if there is no $I_{l}, 0\leq l <n $, in the gap $(I_{k},
I_{n})$ (or $(I_{n}, I_{k})$), where $(I_{k},I_{n})$
denotes the
counter-clockwise gap from $I_{k}$ to $I_{n}$. We
denote them by $I_{L(n)}$ and $I_{R(n)}$.

$I_{n}$ has a successor $I_{n+a}$ if

1. $I_{n-a}$ is a left (or right) predecessor (with $0<a\leq n)$;

2. $f^{a}|_{[I_{n-a},I_{n+a}]}$ (or $f^{a}|_{[I_{n+a},I_{n-a}]}$)
contains no predecessor of $I_{n}$;

3. if $I_{n}$ is to the left (or right) of $I_{n+a}$, then there is no
$I_{k}, 0\leq k< n+a$ in the gap $(I_{n}, I_{n+a})$ (or $(I_{n+a},
I_{n})$).
\end{definition}
Furthermore we define the natural neighborhood $T_{n}$ of $I_{n}$ to be
the biggest closed interval containing $I_{n}$ which contains no
$I_{i},i\in N$, except its nearest predecessor or successor.

Remark: Of course $I_{n}$ can have at most one predecessor on each side.
Moreover $I_{n}$ has at most one successor, denote it by $I_{S(n)}$.
Therefore $T_{n}=[I_{L(n)}, I_{R(n)}]$ if $I_{n}$ has two predecessors
and no successor and $T_{n}=[I_{L(n)}, I_{S(n)}]$ (or
$T_{n}=[I_{S(n)}, I_{R(n)}]$) if $I_{n}$ has a successor.

One can prove the following lemmas ([6], p. 309).
\begin{lemma}
For every $n\in N$, $I_{n}$ can have at most one successor.
\end{lemma}
\begin{lemma}
Assume the interval $I_{n}$ has two predecessors $I_{L(n)}, I_{R(n)}$
and a successor $I_{S(n)}$, If this successor is to the right of
$I_{n}$ then the predecessors of $I_{S(n)}$ are $I_{n}$ and $I_{R(n)}$
and if $I_{S(n)}$ has a successor then this successor must be again to
the right  of $I_{S(n)}$.
\end{lemma}
Remark: This lemma implies that if $I_{n}$ has a successor $I_{S(n)}$
and $I_{S(n)}$ also has a successor $I_{S(S(n))}$ then
$S(n)-n=S(S(n))-S(n)$ and $I_{S(n)}$ is between $I_{n}$ and
$I_{S(S(n))}$. Continuing this if there exists a maximal integer $k$
such
that $I_{S^{i+1}(n)}$ is a successor of $I_{S^{i}(n)}$ for $0\leq i \leq
k-1$, then the intervals $I_{S^{i}(n)}$, $0\leq i \leq k-1$, are
ordered and $f^{a}$ acts as a translation on these intervals, where
$a=S(n)-n$.
\newtheorem{theorem}{Theorem}
\begin{theorem} ([6], p. 310)
Let $n\in N$ and assume that $I_{n}$ has two predecessors $I_{L(n)}$
and $I_{R(n)}$. Let $M_{n}\supset I_{n}$ be an interval contained either
in $[I_{L(n)}, I_{n}]$ or in $[I_{n}, I_{R(n)}]$. Assume that
$\{M_{t_{0}}, M_{t_{0}+1}, ..., M_{n}\}$ are pullbacks of $M_{n}$. If
the intersection multiplicity of this collection is at least $2m$ and
$m\geq 2$ then there exists $t\in \{t_{0}, ..., n\}$ such that

(1) $I_{S(t)}, I_{S^{2}(t)}, ..., I_{S^{2m-2}(t)}$ are defined;

(2) $n=S^{m}(t)$ and $I_{S^{j}(t)}$ is contained in $M_{n}$ for $j=m,
..., 2m-2$.
\end{theorem}
\newtheorem{cor}{Corollary}
\begin{cor}
Assume an interval $T\supset I_{n}$ and $T$ is contained in the
natural neighborhood $T_{n}$ of $I_{n}$. Then the intersection
multiplicity of the pullbacks of $T$ is at most $15$.
\end{cor}
Proof: Consider the pullbacks of $T\cap [I_{L(n)}, I_{n}]$ and $T\cap
[I_{n}, I_{R(n)}]$ seperately. Suppose the intersection multiplicity of
the pullbacks of $T$ is at least $16$. Then either the pullbacks of
$T\cap [I_{L(n)}, I_{n}]$ or the pullbacks of $T\cap [I_{n}, I_{R(n)}]$
has intersection multiplicity $\geq 8$. Take $m=4$, the previous theorem
imples that $I_{S^{2}(n)}$ is contained in $T\cap [I_{n}, I_{R(n)}]$.
This is impossible because of $T\subset T_{n}$.

Now we can prove the following theorem.
\begin{theorem}
Let $f:S^{1}\rightarrow S^{1}$ be an orientation preserving
homeomorphism
with an irrational rotation number. If the logarithm of the cross ratio
distortion under $f$ has bounded variation $B$, then $f$ has no
wandering interval, hence it is topologically conjugate to a rigid
rotation.
\end{theorem}
Proof: Suppose $I$ is a maximal wandering interval for $f$ and
$I_{n}=f^{n}(I), n\geq 0, n\in Z$.
There exists arbitrarily large $n\in N$ and $l, r< n$ such that
$I_{n}\subset (I_{l}, I_{r})$, $I_{k}\cap (I_{l}, I_{r})=\emptyset,
0\leq k<n$, and $|I_{n}|\leq \min \{|I_{l}|, |I_{r}|\}$. This property
is proved as follows. Pick up $I_{l}$
and $I_{r}$ such that the gap $(I_{l}, I_{r})$ contains no $I_{k}$ for
$0\leq k\leq \max \{l, r\}$. By the density of any
orbit under an irrational rotation, there exists $I_{n}$ first gets into
the gap $(I_{l}, I_{r})$. If $|I_{n}|\leq \min \{|I_{l}|, |I_{r}|\}$
then it is done, otherwise replace $I_{r}$ by $I_{n}$ and go on.
Since the sum of the lenths of $I_{k}$ is bounded, eventually we will
get $|I_{n}|\leq |I_{r}|$, furthermore we get $|I_{n}|\leq |I_{l}|$.
We have seen $I_{l}$ and $I_{r}$ are two predecessors of $I_{n}$.

Let $T_{n}$ be the natural neighborhood of $I_{n}$. If $I_{n}$ has no
successor, then $T_{n}=[I_{l}, I_{r}]$. By the corollary 1, the
intersection multiplicity of the pullbacks of $T_{n}$ is bounded by
$15$. Use the Macroscopic Koebe Distortion Principle, we get
a bigger wandering interval $J$ strickly containing $I$. This
contradicts with the maximality of $I$. Hence $I_{n}$ has a
successor $I_{s(n)}$. Use the same way as the above, we get
$|I_{s(n)}|<|I_{n}|$. Inductively we get infinitely many
successors $I_{s^{i}(n)}, i=1,2,...$, and by theorem 1, all successors
are contained in $[I_{l}, I_{r}]$ and are ordered. Moreover
$s^{i}(n)-s^{i-1}(n)$ is a constant $a=s(n)-n$. It follows
$I_{s^{i}(n)}$ converges to a fixed point of $f^{a}$ as $i\rightarrow
\infty $. This contradicts with $f$ has no periodic points.

\vspace{.2in}
${\bf Proof\; of\; Theorem \;C:}$ Let $I$ be a maximal wandering
interval for
$f$ and $I_{n}=f^{n}(I), n\geq 0, n\in Z$. Let $T_{n}$ be the natural
neighborhood of $I_{n}$. The intersection multiplicity of the pullbacks
of $T_{n}$ is bounded by $15$. By the Prop. 4 the cross ratio distortion
of $f^{n}$ on $T_{0}$ is uniformly bounded by a constant $B$. The rest
of the proof follows the proof of the above theorem.

\vspace{.2in}
The remainder of this section explains why the conditions of
the
theorem C is weaker than Denjoy's condition and [5]'s condition.
It is almost trivial that Denjoy's condition implies the conditions of
the theorem C.
\begin{prop}
Let $h: I\rightarrow R^{1}$ be a $C^{1}$ smooth function and $logh'$ is
of bounded variation, then $logh'$ is of bounded Zygmund variation and
bounded quadratic variation.
\end{prop}
Proof: By the triangle inequality, the Zygmund variation of $logh'$ is
no more than the variation of $logh'$ on the interval $I$.
Let $M$ be the maximal value of $|logh'|$ on the interval $I$,
then the quadratic variation of $logh'$ is no more than $2M$
multiplied by the variation of $logh'$ on the interval $I$.

Clearly the Zygmund condition implies bounded Zygmund variation.
Furthermore, the Zygmund condition implies $\alpha $-H\"older continuous
for $0<\alpha <1$. The $1/2$-H\"older continuity implies the bounded
quadratic variation.
\begin{lemma}
If $\phi :I\rightarrow R^{1}$ satisfies the Zygmund condition:
there exists $B>0$ such that
\[\sup_{x,t}|\frac{\phi (x+t)+\phi (x-t)-2\phi (x)}{t}|\leq B,\]
then $\phi $ is $\alpha $-h\"older continuous for any $0<\alpha <1$.
\end{lemma}
Proof: Denote $D(x,t)=\frac{\phi (x+t)-\phi (x)}{t}$. Then
\[D(x,t/2)+D(x+t/2,t/2)=2D(x,t),\;|D(x,t/2)-D(x+t/2,t/2)|\leq B,\]
\[D(x,t/4)+D(x+t/4,t/4)=2D(x,t/2),\;|D(x,t/4)-D(x+t/4,t/4)|\leq B,\]
\[\cdot \]
\[\cdot \]
\[\cdot \]
\[D(x,\frac{t}{2^{n}})+D(x+\frac{t}{2^{n}},\frac{t}{2^{n}})
=2D(x,\frac{t}{2^{n-1}}),\;
 |D(x,\frac{t}{2^{n}})-D(x+\frac{t}{2^{n}},\frac{t}{2^{n}})|\leq B.\]
These give us
\[|D(x,t/2^{n})|\leq |D(x,t)|+nB, \] i. e.,
\[|\frac{\phi (x+t/2^{n})-\phi (x)}{t/2^{n}}|\leq |D(x,t)|+nB.\]
Then \[|\frac{\phi (x+t/2^{n})-\phi (x)}{(t/2^{n})^{\alpha }}|\leq
(|D(x,t)|/n+B)n(|t|/2^{n})^{1-\alpha },\]
which tells us that $\phi $ is $\alpha $-H\"older continuous for any
$0<\alpha <1$.
\begin{prop}
If $\phi :I\rightarrow R^{1}$ satisfies the Zygmund condition, then
$\phi $ is of bounded Zygmund variation and bounded quadratic variation
over the interval $I$.
\end{prop}
\section{Three examples}
In the introduction it is mentioned that there exists an example
satisfying
Denjoy's bounded variation condition but [5]'s Zygmund condition
and vice versa. In this section, we will give these two examples and
also we will give an example to show that there is an example being
of bounded quadratic variation but not being of bounded Zygmund
variation.
\newtheorem{example}{Example}
\begin{example} Let $\phi :[-1,\;1]\rightarrow [-1,\;1]$ be the
following function \[\phi (x)=x,\;\;\;x\in [-1,\;0],\]
	 \[\phi (x)=\sqrt {x},\;\;\;x\in (0,\;1].\]
Clearly $\phi $ is monotone hence it is of bounded variation, but the
Zygmund condition fails since the right derivative of $\phi $ at the
point $0$ is infinite but the left derivative is $1$.
\end{example}
\begin{example}Let $\phi _{0}(x)=2x$ for $x\in [0,
\frac {1}{2}]$ and $\phi _{0}(x)=2-2x$ for $x\in [\frac {1}{2}, 1]$. And
let \[\phi _{n}=\frac {\phi (2^{n}x-i)}{2^{n}} \;for \; x\in [\frac
{i}{2^{n}}, \frac{i+1}{2^{n}}], \]
where $i=0, 1, ..., 2^{n}-1$.

Let $\phi (x)=\sum _{n=0}^{\infty } \phi _{n} (x)$.
$\phi $ is differentiable at a set of measure $0$ only. It can't be
of bounded variation, otherwise it is differentiable almost
everywhere which is a contradiction. [7] and [8] study the general
theory about the differentiability of a function satisfying Zygmund
condition.
\end{example}
\begin{example}Let $\phi : [0,\;1]\rightarrow [0,\;1]$ is defined by
the figure 1.
\begin{figure}[t]
\par \centerline{ \hbox{\psfig{figure=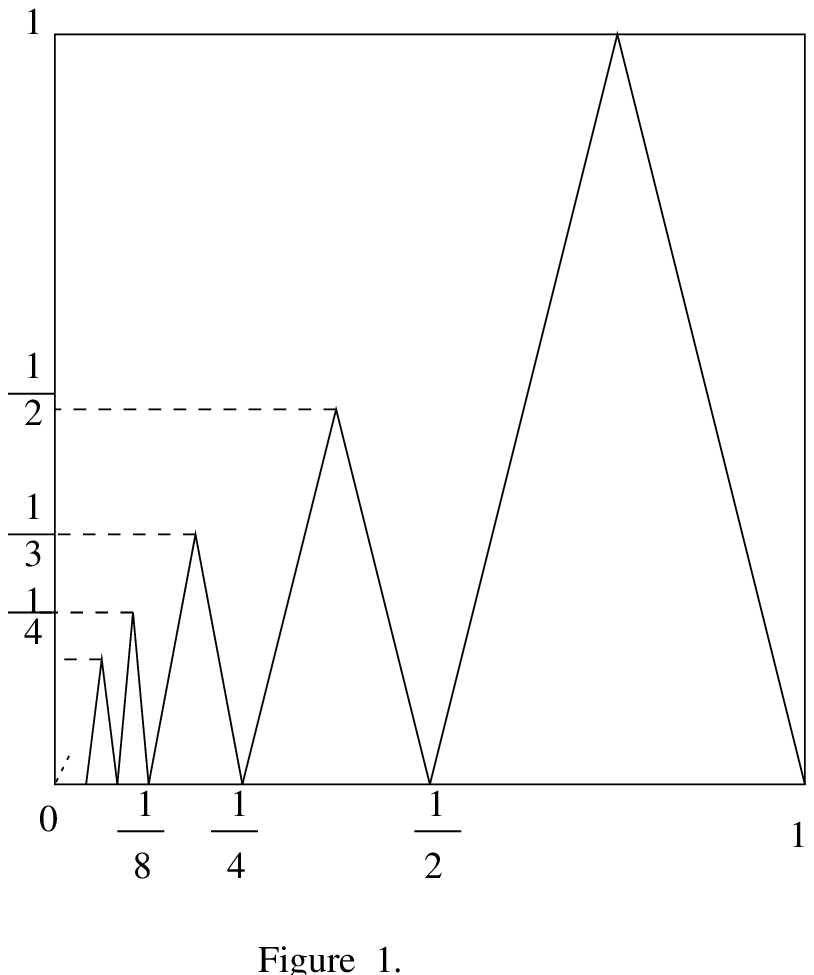,width=3.0in }}}\par
\end{figure}
It is easy to get the quadratic variation of $\phi $ on $[0,\;1]$ is
equal to $\sum_{n=1}^{\infty }\frac{2}{n^{2}}$ which is finite. But the
difference between the left derivative and the right derivative
of $\phi $ at $1/2^{n}$ is equal to
\[\frac{2^{n+2}}{n+1}-\frac{2^{n+1}}{n}=\frac{2^{n+1}}{n}\frac{n-1}{n+1}\]
which tends to $\infty $ as $n\rightarrow \infty $. Therefore it has
bounded quadratic variation but has no bounded Zygmund variation.
\end{example}
Actually the left question is to study whether or not the
bounded Zygmund variation property implies the bounded quadratic
variation property.
\section{Appendix}
The nonexistence of wandering domains for any rational map of the
complex sphere was proved by Dennis Sullivan in 1985 [9]. The
analogue of this theorem for one-dimensional dynamical systems was done
for certain smooth multimodal maps by Martens, de Melo and van Strien in
1992 [10]. In the latest publication [6], the smooth condition used by
de Melo and van Strien is that a multimodal map piecewise satisfies
$C^{1+b.v}$
(or $C^{1+Z}$) and the map can be written as a power map $(x\mapsto
|x|^{\alpha }, \alpha >1)$ composed by a $C^{1+b.v}$ (or $C^{1+Z}$)
diffeomorphism around every turning point. Combine the analysis work
of getting a bound of cross ratio distortions in this paper and
the combinatorial machinery on wandering intervals in [10] or
[6, p. 308-312], we can get a weak version of Martens, de Melo
and van Strien's theorem of no wandering intervals for multimaodal maps.
Before we state the theorem, let us give the definition of a wandering
interval for a multimodal map of an interval.

\begin{definition} Let $f: I\rightarrow I$ be a continuous map of an
interval $I$.
An open interval $J\subset I$ is called a wandering interval of $f$ if

1) $f^{n}(J)\cap f^{m}(J) =\emptyset$ for any $n\neq m, n,m\in N$;

2) $f^{n}(J)$ does not converge to a periodic orbit.
\end{definition}
\begin{theorem}$[11]$ Let $f: I\rightarrow I$ be a $C^{1}$
smooth map satisfying

1) $f$ is $C^{1+b.Z.v+b.q.v}$ away from critical points;

2) Let $K_{f}$ be the set of critical points of $f$. For each $x_{0}\in
K_{f}$, there exist $\alpha >1$, a neighborhood $U(x_{0})$ of $x_{0}$
and a $C^{1+b.Z.v+b.q.v}$ diffeomorphism $\phi :U(X_{0})\rightarrow
(-1, 1)$ such that $\phi (x_{0})=0$ and \[f(x)=f(x_{0})\pm |\phi
(x)|^{\alpha }, \forall x \in U(x_{0}).\]
Then $f$ has no wandering intervals.
\end{theorem}
Norton, Sullivan and Velling [12, 13 and 14] have begun the work of
generalizing
the setting of Denjoy's theorem to two dimensional dynamical systems by
considering diffeomorphisms of the torus. The quasiconformal theory
has found a place there.

\vspace{.2in}
\noindent{\em Acknowledgements}.
Both authors wish to thank Frederick P. Gardiner for his help in
writing, and thank Yunping Jiang and Meiyu Su for some discussions.
We are grateful to John Milnor for his suggestions to improve the writing 
and especially for his simple setting of the proof of Prop. 2.

\end{document}

%% file: imsmark.tex
\def\IMSmarkvadjust{0 pt}
\def\IMSmarkhadjust{0 pt}
\def\IMSmarkhpadding{0 pt}
\def\IMSpubltext{Published in modified form:}
\def\SBIMSMark#1#2#3{
 \font\SBF=cmss10 at 10 true pt
 \font\SBI=cmssi10 at 10 true pt
 \setbox0=\hbox{\SBF \hbox to \IMSmarkhpadding{\relax}
                Stony Brook IMS Preprint \##1}
 \setbox2=\hbox to \wd0{\hfil \SBI #2}
 \setbox4=\hbox to \wd0{\hfil \SBI #3}
 \setbox6=\hbox to \wd0{\hss
             \vbox{\hsize=\wd0 \parskip=0pt \baselineskip=10 true pt
                   \copy0 \break%
                   \copy2 \break%
                   \copy4 \break}}
 \dimen0=\ht6   \advance\dimen0 by \vsize \advance\dimen0 by 8 true pt
                \advance\dimen0 by -\pagetotal
	        \advance\dimen0 by \IMSmarkvadjust
 \dimen2=\hsize \advance\dimen2 by .25 true in
	        \advance\dimen2 by \IMSmarkhadjust

  \openin2=publishd.tex
  \ifeof2\setbox0=\hbox to 0pt{}
  \else 
     \setbox0=\hbox to 3.1 true in{
                \vbox to \ht6{\hsize=3 true in \parskip=0pt  \noindent  
                {\SBI \IMSpubltext}\hfil\break
                \input publishd.tex 
                \vfill}}
  \fi
  \closein2
  \ht0=0pt \dp0=0pt
 \ht6=0pt \dp6=0pt
 \setbox8=\vbox to \dimen0{\vfill \hbox to \dimen2{\copy0 \hss \copy6}}
 \ht8=0pt \dp8=0pt \wd8=0pt
 \copy8
 \message{*** Stony Brook IMS Preprint #1, #2. #3 ***}
}

%% file: publishd.tex
{\it Ergod. Th. \& Dynam. Sys.}~{\bf 17} (1997), 173--186